\documentclass[graybox]{svmult}

\usepackage{type1cm}        
\usepackage{makeidx}         
\usepackage{graphicx}        
\usepackage{multicol}        
\usepackage[bottom]{footmisc}
\usepackage{multirow}

\usepackage{newtxtext}       %
\usepackage{newtxmath}       

\usepackage{color}
\usepackage{ulem}
\usepackage{tikz}
\usepackage{bbold}
\usepackage{enumitem}

\makeindex             

\usepackage[capitalise]{cleveref}

\begin{document}

\title*{Nonlinear Monolithic Two-Level Schwarz Methods for the Navier-Stokes Equations}
\titlerunning{Nonlinear Two-Level Schwarz Methods for Navier-Stokes}
\author{Axel Klawonn\orcidID{0000-0003-4765-7387} \\ Martin Lanser\orcidID{0000-0002-4232-9395}}
\institute{Axel Klawonn$^{1,2}$, Martin Lanser$^{1,2}$ \at $^1$Department of Mathematics and Computer Science, Division of Mathematics, University of Cologne, Weyertal 86-90, 50931 Cologne, Germany, \email{axel.klawonn@uni-koeln.de, martin.lanser@uni-koeln.de}, {url:~\url{https://www.numerik.uni-koeln.de}} \at $^2$Center for Data and Simulation Science, University of Cologne, Germany, {url:~\url{https://www.cds.uni-koeln.de}}}
%
%
\maketitle

\abstract{Nonlinear domain decomposition methods became popular in recent years since they can improve the nonlinear convergence behavior of Newton's method significantly for many complex problems. In this article, a nonlinear two-level Schwarz approach is considered and, for the first time, equipped with monolithic GDSW (Generalized Dryja-Smith-Widlund) coarse basis functions for the Navier-Stokes equations. Results for lid-driven cavity problems with high Reynolds numbers are presented and compared with classical global Newton's method equipped with a linear Schwarz preconditioner. Different options, for example, local pressure corrections on the subdomain and recycling of coarse basis functions are discussed in the nonlinear Schwarz approach for the first time.}

\section{Introduction and Nonlinear Problems}
\label{sec:intro}

In recent years, nonlinear preconditioning techniques became popular to improve the convergence speed of nonlinear solvers for complex partial differential equations. Many of these nonlinear preconditioners are based on the ideas and basic principles of linear domain decomposition methods (DDMs) and consequently they are denoted as nonlinear DDMs. In nonlinear DDMs, e.g.,~\cite{schwarz1,schwarz2,schwarz5,nlfeti1,nlfeti2,nlfeti3,nlfeti4,nlfeti5,raspen,mspin1,mspin2,krause,schwarz6}, the discretized nonlinear problem is restricted to smaller local ones on subdomains and the local solutions of those are recombined to the original one in a Newton-like iteration, which hopefully converges faster than classical Newton's method applied to the original problem. As in the linear case it often is beneficial to include a (nonlinear) coarse space for robustness and convergence improvement, which leads to two-level nonlinear DDMs. In this article, we consider some variants of the nonlinear two-level Schwarz approaches first introduced in~\cite{schwarz1}; older two-level nonlinear Schwarz approaches for less general coarse spaces can be found in~\cite{schwarz6,raspen}. Here, for the first time, we combine two-level nonlinear Schwarz methods with a monolithic RGDSW (Reduced Generalized Dryja-Smith-Widlund) coarse space for Navier-Stokes problems (see~\cite{schwarz2,schwarz3}) and apply it to the well-known lid-driven cavity flow test problem considering high Reynolds numbers. 
We compare our solvers with classical Newton-Krylov-Schwarz approaches equipped with a simple globalization, namely a simple backtracking approach.   

\paragraph{\bf Lid-driven Cavity Problem}
As stated above, we consider the classical lid-driven cavity flow problem modeled by the stationary, dimensionless, and incompressible Navier-Stokes equations in two-dimensions
\begin{equation*}
\begin{array}{rll}
	-\frac{1}{Re} \Delta v + (v \cdot \nabla)v + \nabla p &= 0 \; &{\rm in}\; \Omega\\
{\rm div}(v) &= 0 \; &{\rm in} \; \Omega \\
v &= v_0 \; &{\rm on} \; \partial \Omega \\
p &=0 \; &{\rm in} \; (x,y)=(0,0). 
\end{array}	
\end{equation*}
Here, $v$ denotes the velocity, $p$ the pressure, and $Re$ the Reynolds number. In case of the lid-driven cavity problem $\Omega$ is the unity square and $v_0=[1,0]$ for the Dirichlet boundary constraints on the lid, that is, the upper boundary of $\Omega$ where $y$ is equal to one. On all other parts of the boundary $\partial \Omega$ we have $v_0=[0,0]$. We usually use $v^{(0)}|_{\partial \Omega}=v_0$ and $v^{(0)}=0$ within $\Omega$ and $p^{(0)}=0$ as initial value for Newton's method, regardless if we use a classical Newton-Krylov-Schwarz approach or a modern nonlinear Schwarz method. Therefore, since $v^{(0)}$ already fulfills all boundary constraints, we can enforce a zero Dirichlet boundary condition on $\partial \Omega$ for the velocity of the Newton update in each Newton iteration. Note that for the pressure we enforce a Dirichlet boundary condition in a single corner, that is, in $(x,y) = (0,0)$ to obtain a unique solution. 

\section{Nonlinear Two-Level Schwarz Methods}
\label{sec:nl}

Let us briefly describe the two-level nonlinear Schwarz methods considered in our experiments, which were first introduced in~\cite{schwarz1} and extended in~\cite{schwarz4} for more general coarse spaces. We assume to have a nonlinear problem 
\begin{equation}
	F(u)=0
\label{eq:basic}
\end{equation}
which originates from a finite element discretization of a given nonlinear partial differential equation. In the following, we define nonlinearly left-preconditioned systems of the form
\begin{equation}
	\mathcal{F}(u) := G(F(u)) = 0,
	\label{eq:nl}
\end{equation}
where the left preconditioner $G$ is given implicitly by a domain decomposition approach. The nonlinearly preconditioned problem is then solved with Newton's method instead of solving the original formulation from~\cref{eq:basic} directly. The goal of this approach is to accelerate the nonlinear convergence of Newton's method by a proper choice of $G$. 

We first consider a decomposition of the computational domain $\Omega \subset \mathbb{R}^d,\; d=2,3$, into $N$ nonoverlapping subdomains $\Omega_i',\; i=1,...,N$, such that $\overline{\Omega} = \bigcup\limits_{i=1}^N \overline{\Omega}_i'$.
By adding $m$ rows of finite elements around the boundary of each subdomain, we obtain an overlapping domain decomposition $\overline{\Omega} = \bigcup\limits_{i=1}^N \overline{\Omega_i}$ with overlap $\delta =m \cdot h$, where $h$ is the typical diameter of a finite element. 

For the first level of the nonlinear Schwarz preconditioner, we define local nonlinear corrections $T_i(u)$ by solving 
\begin{equation}
	R_i F(u- P_i T_i(u))=0, \; i=1,...,N,
	\label{eq:local_corr}
\end{equation}
where $R_i: V \rightarrow V_i$ is the restriction from the global finite element space $V$ to the local space $V_i$ belonging to the overlapping subdomain $\Omega_i, \; i=1,...,N$. Here, $P_i: V_i \rightarrow V, \; i=1,...,N$, is a corresponding prolongation operator; we always use the symmetric choice of $P_i := R_i^T$ throughout this article. 

For the definition of the second and coarse level, we assume to have a coarse space $V_0$ and a prolongation operator $P_0: V_0 \rightarrow V$. In practice, we usually compute the operator $P_0$, which is a column-wise collection of coarse basis functions discretized in $V$ and afterwards define 
$$
V_0 := \{ v_0 = P_0^T v \; | \; v \in V \}.
$$
Let us note that we are usually interested in coarse spaces which can be built without having a coarse discretization of $\Omega$. Nonetheless, if such coarse finite elements defining $V_0$, the corresponding operators $R_0$ and $P_0$ can be built by simply interpolating the finite element shape functions from $V_0$ to the fine space $V$. Details on the construction of $P_0$ for the Navier-Stokes equations without having a coarse discretization are given in \cref{sec:coarse}. With $R_0 := P_0^T$, we can define a nonlinear coarse correction $T_0(u)$ by
\begin{equation}
	R_0 F(u- P_0 T_0(u))=0.
	\label{eq:coarse_corr}
\end{equation}
Now we have all ingredients to define new left-preconditioned nonlinear operators and we start with an additive variant. To obtain an additive nonlinear two-level Schwarz method, the nonlinear problem
\begin{equation}
\mathcal{F}_a(u) := \sum\limits_{i=1}^N P_i T_i(u) + P_0 T_0(u)
\label{eq:1}
\end{equation}
is solved using Newton's method. We refer to this method as two-level ASPEN (Additive Schwarz Preconditioned Exact Newton) based on the one-level approaches from~\cite{schwarz5}, or, even simpler, as additive nonlinear two-level Schwarz method. Similarly, a hybrid nonlinear two-level Schwarz method can be obtained by solving
\begin{equation}
\mathcal{F}_h(u) := \sum\limits_{i=1}^N P_i T_i(u-P_0 T_0(u)) + P_0 T_0(u)
\label{eq:2}
\end{equation}
with Newton's method. Here, the coupling between the levels is multiplicative while the subdomains are still coupled additively against each other. For more details and further nonlinear two-level variants, we refer to~\cite{schwarz1,schwarz4,raspen}.

Let us remark that in each Newton iteration a linear system with the Jacobian $D\mathcal{F}_X(u^{(k)})$ has to be solved iteratively with, for example, GMRES (Generalized Minimal Residual), where $X \in \{a,h \}$. Fortunately, $D\mathcal{F}_X(u^{(k)})$ has already the favorable structure of a linear two-level Schwarz preconditioner, either additive or hybrid, times the Jacobian $DF(\cdot)$. Therefore, no additional linear preconditioner is necessary to solve the linearized systems efficiently. For details on the exact shape and the different linearization points of the local and the coarse part of the preconditioner, we refer to~\cite{schwarz1}. For experts, let us remark that we use always the exakt Jacobian and not the approximation used, for example, in the well-known ASPIN approach~\cite{schwarz5}. All in all, nonlinear two-level Schwarz algorithms have a special structure. First, there is an outer Newton loop where in each iteration a linearized system with the matrix $D \mathcal{F}_X(u^{(k)})$ is solved iteratively using GMRES. Second, to compute the residual $\mathcal{F}_X(u^{(k)})$ and the Jacobian $D \mathcal{F}_X(u^{(k)})$, all corrections $T_i(u^{(k)}),\; i=0,...,N$, have to be computed in each outer Newton iteration, which is done by solving \cref{eq:local_corr} and \cref{eq:coarse_corr} with Newton's method. This immediately leads to inner Newton iterations which are denoted by local inner loops in case of the local corrections and the coarse inner loop in case of the coarse correction. Let us note that all local inner loops can be easily run in parallel but the coarse inner loop can only be carried out in parallel in the case of an additive coupling between levels. In the hybrid variant, the coarse loop is carried out before the subdomain corrections are computed. We provide a brief algorithmic overview of two-level nonlinear Schwarz methods in \cref{alg:1}. 

\begin{figure}
\begin{itemize}
\item[] {\bf Outer Newton loop} over $k=0,1,2,...$\\[-4ex]
\begin{itemize}[itemsep=2pt]
\item[] {\bf Inner coarse loop} to compute $T_0(u^{(k)})$ (solving \cref{eq:coarse_corr}  with Newton's method)
\item[] {\bf Inner local loops} to compute $T_i(u^{(k)})$ {\bf OR} $T_i(u^{(k)}-P_0T_0(u^{(k)}))$, $i=1,...,N$\\[-3.5ex]
\item[] (solving \cref{eq:local_corr}  with Newton's method; parallelizable)
\item[] {\bf Set Up} $\mathcal{F}_X(u^{(k)})$ and $D\mathcal{F}_X(u^{(k)})$ with $X \in \{a,h\}$ using $T_i(u^{(k)}), \; i=0,...,N$
\item[] {\bf Iterative solution with GMRES} of $D\mathcal{F}_X(u^{(k)}) \, \delta u^{(k)} = \mathcal{F}_X(u^{(k)})$
\item[] {\bf Update} $u^{(k+1)} = u^{(k)} - \delta u^{(k)}$\\[-4ex]
\end{itemize}
\item[] {\bf End} of outer loop	
\end{itemize}
\caption{Brief algorithmic overview of additive or hybrid nonlinear two-level Schwarz.}
\label{alg:1}	
\end{figure}

To get a good estimate of the performance of our nonlinear Schwarz methods, we compare against classical Newton-Krylov-Schwarz, where Newton's method is applied to solve $F(u)=0$ directly and all linear systems are solved using GMRES and a linear two-level Schwarz preconditioner, either hybrid or additive. %
In our comparison, we always use the same coarse spaces for Newton-Krylov-Schwarz and nonlinear Schwarz. %

\section{Coarse Spaces}
\label{sec:coarse}

For our Navier-Stokes problems we use the monolithic basis functions introduced in~\cite{schwarz2}, which are based on the GDSW coarse spaces; see~\cite{gdsw1} for scalar elliptic problems. 
We only provide a brief description here. Let us first assume that we have a symmetric and positive definite matrix $K$ which is obtained from the discretization of a scalar linear partial differential equation, for example, a linear diffusion equation with sufficient boundary constraints of Dirichlet type. Thus, we have a single degree of freedom in each node of the finite element mesh and, for the moment, do not have to distinguish between nodes and degrees of freedom. We partition $K$ into interface and interior nodes
$$
K = \begin{bmatrix}
	K_{II} & K_{I \Gamma}\\
	K_{\Gamma I} & K_{\Gamma \Gamma}
\end{bmatrix},
$$
where the interface $\Gamma$ contains all nodes on the boundaries of the nonoverlapping subdomains $\Omega_i',\; i=1,...,N,$ except the nodes belonging to a part of $\partial \Omega$ where a Dirichlet boundary condition is imposed. The set $I$ contains all remaining nodes. Furthermore, we partition the interface $\Gamma$ into $n_c$ overlapping or nonoverlapping patches $\Gamma_j,\; j=1,...,n_c$, and define associated functions or, more precisely, vectors $\Phi_\Gamma^j, \; j=1,...,n_c$, which are defined on $\Gamma$ but are only nonzero on the corresponding patch $\Gamma_j,\; j=1,...,n_c$. Additionally, we construct $\Phi_\Gamma^j$ such that we obtain a partition of unity $\sum_{j=1}^{n_c} \Phi_\Gamma^j = 1$ on the complete interface. %
Using the discrete energy minimizing extension
$$
\Phi_I^j := - K_{II}^{-1} K_{I \Gamma} \Phi_\Gamma^j
$$
we define the coarse basis function belonging to the patch $\Gamma_j, \; j=1,...,n_c$, by $$\Phi^j := \begin{pmatrix}
\Phi_I^j \\ \Phi_\Gamma^j	
\end{pmatrix}\quad  {\rm and \; by} \quad 
\Phi := \begin{bmatrix}
 	\Phi^1,...,\Phi^{n_c}
 \end{bmatrix}
$$
the matrix which contains all coarse basis functions. Let us remark that $P_0 := \Phi$ gives us the prolongation from the coarse space $V_0$ to the fine space $V$ and that we have a partition of unity
$$
\sum_{j=1}^{n_c} \Phi^j = \mathbb{1}
$$
on $\Omega$ with $\mathbb{1}$ being a vector of ones, since $K$ is symmetric positive definite. 

Using this concept, different partitions into patches $\Gamma_j$ and different definitions of $\Phi_\Gamma^j$ are possible. Partitioning $\Gamma$ into $n_e$ edges and $n_v$ vertices in two dimensions and defining $\Phi_\Gamma^j, \; j=1,...,n_c,$ with $n_c=n_v+n_e$ to be one on the corresponding edge or vertex leads to the classical GDSW coarse space. Alternatively, we can define one patch $\Gamma_j, \; j=1,...,n_v$, for each vertex, consisting of the vertex itself and all adjacent edges. Then, $\Phi_\Gamma^j$ can be defined to be one in the vertex and to linearly decrease to zero along the adjacent edges. This gives a smaller coarse space of size $n_v$ but also fulfills all necessary properties. This coarse space is called MsFEM-D coarse space in \cite{schwarz1} and can also be interpreted as one instance of the class of reduced GDSW (RGDSW) coarse spaces; see~\cite{schwarz3,rgdsw} for more details on RGDSW. In the present article, we denote this coarse space as RGDSW coarse space of type A. An alternative approach to define $\Phi_\Gamma^j$ is to set it to $1$ in the vertex and to $0.5$ on the adjacent edges. Let us note that if an edge ends at the boundary $\partial \Omega$ we set it to 1 instead of $0.5$. This coarse space can be implemented in an algebraic way without information about the geometry of the subdomains. We refer to this as RGDSW type B coarse space.

In general, for a nonlinear problem, one can use the Jacobian $K := DF(u^{(k)})$ in the $k$-th linearization $u^{(k)}$ of Newton's method to compute the coarse basis functions. We consider two approaches: 1) compute $\Phi$ once using $K:= DF(u^{(0)})$ in the initial value $u^{(0)}$ and then recycle $\Phi$ in all further steps or 2) recompute $\Phi$ with $K := DF(u^{(k)})$ at the beginning of each Newton step. Both versions can be used either for Newton-Krylov-Schwarz or nonlinear two-level Schwarz, where $\Phi$ can be recomputed at the beginning of each outer iteration. 

Considering the Navier-Stokes equations in two dimensions, we have three degrees of freedom for each node, that is, two velocities $v=(v_x,v_y)$ and the pressure $p$. To build a partition of unity on the interface, we now define three basis functions $\Phi_\Gamma^{j,i},\;i=1,2,3,\; j=1,...,n_c,$ for each patch $\Gamma_j, \; j=1,...,n_c$; one for each degree of freedom (velocities and pressure). Sorting the degrees of freedom appropriately, that is, first all $x$-velocities, then all $y$-velocities, and finally all pressures, we simply have
$$
\Phi_\Gamma^{j,1}=\begin{bmatrix}
\Phi_\Gamma^{j} \\ 0 \\ 0	
\end{bmatrix}, \quad
\Phi_\Gamma^{j,2}=\begin{bmatrix}
0\\ \Phi_\Gamma^{j} \\ 0	
\end{bmatrix}, \quad {\rm and} \quad
\Phi_\Gamma^{j,3}=\begin{bmatrix}
 0 \\ 0 \\ \Phi_\Gamma^{j}	
\end{bmatrix}
$$
and, of course,
$$
\sum\limits_{i=1}^3 \sum\limits_{j=1}^{n_c} \Phi_\Gamma^{j,i} =\mathbb{1}
$$
on the interface. For the extension to the interior degrees of freedom, the complete block matrix
$$
K = DF(u^{(k)}) = \begin{bmatrix}
 A(u^{(k)}) & B^T\\ B & 0	
 \end{bmatrix}
$$ 
from the linearization of the discretized Navier-Stokes equation is used, which is clearly a monolithic approach to compute the coarse basis. Here, as usual, only the velocity part $A$ depends on the current velocities and pressures $u^{(k)} = \left(v^{(k)},p^{(k)}\right)$. We note that, since $K$ is no longer a symmetric positive definite matrix, we have a lack of theory and cannot show that the coarse basis functions $\Phi$ build a partition of unity on the interior degrees of freedom. Nonetheless, in case of Newton-Krylov-Schwarz, this approach has already shown good performance; see~\cite{schwarz2}. Let us remark that the resulting $\Phi$ can also be written as
$$
\Phi = \begin{bmatrix}
\Phi_{vv} & \Phi_{vp}\\
\Phi_{pv} & \Phi_{pp}	
\end{bmatrix},
$$
where the first column contains all basis functions belonging to the velocities $v$ and the second column contains all basis functions belonging to the pressure variables. In our computations, we use a modification by deleting the coupling blocks, that is,
$$
\widetilde{\Phi} = \begin{bmatrix}
\Phi_{vv} & 0\\
0 & \Phi_{pp}	
\end{bmatrix}.
$$
These modified coarse basis functions yield better results in the case of nonlinear Schwarz methods and have already been used in~\cite{schwarz2}. For the definition of the patches $\Gamma_j$ and corresponding scalar basis functions $\Phi_\Gamma^j$ we use both RGDSW type approaches suggested above.

To further improve the linear convergence, we additionally enforce a zero pressure average for the subdomain corrections in the linear solves. Following~\cite{schwarz2}, we define the local projection $$\overline{P}_i := I_p - a_i^T (a_i \, a_i^T)^{-1} a_i,\; i=1,...,N,$$ on the pressure variables of the subdomain $\Omega_i,\; i=1,...,N$, with $a_i=1$ is a vector of ones; one for each local pressure variable. We then define the local projection
$$
\overline{\mathcal{P}}_i := \begin{bmatrix}
I_v	& 0\\ 0 & \overline{P}_i
\end{bmatrix}
$$
for each subdomain $\Omega_i,\;i=1,...,N$, which is the identity for all local velocities but enforces a zero average in the pressure.

Usually, we prefer to use restricted versions of nonlinear Schwarz. In the restricted variants, the local projections $P_i,\; i=1,...,N,$ which are used to recombine the nonlinear local corrections are replaced by $\widetilde{P}_i,\; i=1,...,N$, where the $\widetilde{P}_i$ fulfill a partition of unity property on $\Omega$, that is, $\sum\limits_{i=1}^N \widetilde{P}_i R_i = I$.
 In our numerical results, due to space limitatons, we only consider the hybrid two-level restricted nonlinear Schwarz approach with or without local pressure projections, that is, we either solve
\begin{equation}
\overline{\mathcal{F}}_h(u) := \sum\limits_{i=1}^N \widetilde{P}_i \overline{\mathcal{P}}_i T_i(u-P_0 T_0(u)) + P_0 T_0(u) = 0
\label{eq:3}
\end{equation}
or 
\begin{equation}
\mathcal{F}_h(u) := \sum\limits_{i=1}^N \widetilde{P}_i T_i(u-P_0 T_0(u)) + P_0 T_0(u) = 0
\label{eq:4}
\end{equation}
with Newton's method. Let us remark that the local projections have successfully been used to reduce the number of GMRES iterations in the case of Newton-Krylov-Schwarz but now, in the case of nonlinear Schwarz, $\overline{\mathcal{P}}_i$ actually appears in the nonlinear formulation and might influence the convergence of Newton's method.  

Let us give some details on the initial value of Newton's method. We decided to use an initial value which fulfills all Dirichlet boundary constraints such that we can use homogeneous boundary conditions in all linear solves. We thus choose $v^{(0)}|_{\partial \Omega}=v_0$ and $v^{(0)}=0$ elsewhere. We use $p^{(0)}=0$ for the pressure. Since we have zero Dirichlet boundary conditions for the velocity now, we only choose coarse basis functions in the vertices which are not part of the boundary for both, velocity as well as for the pressure.

To summarize, we vary three different aspects in the numerical tests: 1) the coarse space using two different types of RGDSW, 2) using or not using local pressure projections $\overline{P}_i$, and 3) recycling the coarse basis functions $\widetilde{\Phi}$ or recomputing them in each Newton iteration.

\section{Numerical Results for the Lid-driven Cavity Flow Problem}
\label{sec:results} 
As already stated, we consider the classical lid-driven cavity problem; see \cref{sec:intro} for details. As a classical approach, we use Newton-Krylov-Schwarz with a linear hybrid restricted two-level Schwarz preconditioner. We always use the monolithic RGDSW basis functions described above to build the second level using either type A or B to define them on the interface. The linear preconditioner is not varied since it will not affect the nonlinear convergence in case of Newton-Krylov-Schwarz. We choose only the best possible option which is, in this specific case,  recycling the coarse basis functions $\widetilde{\Phi}$ and using the local pressure projections $\overline{\mathcal{P}}_i,\; i=1,...,N$. As a nonlinear Schwarz approach, we use the hybrid restricted two-level Schwarz approach with or without local pressure projection, that is, we  consider one of the residual formulations~\cref{eq:3} or~\cref{eq:4} and solve with Newton's method. We can optionally also use recycling of the coarse basis functions in nonlinear Schwarz methods and also compare both initial values.

We stop the nonlinear iteration when the relative residual $||F(u^{(k)})||_2/||F(u^{(0)})||_2$ is smaller than $10^{-6}$. All inner local and coarse nonlinear iterations in nonlinear Schwarz are stopped early when a reduction of the relative residual of at least three orders of magnitude is reached. Furthermore, all iterative solves with GMRES are stopped when the relative residual is below $10^{-8}$. We test for different Reynolds numbers $Re$ up to $2500$.

\begin{figure}
\centering
\includegraphics[width=0.9\textwidth]{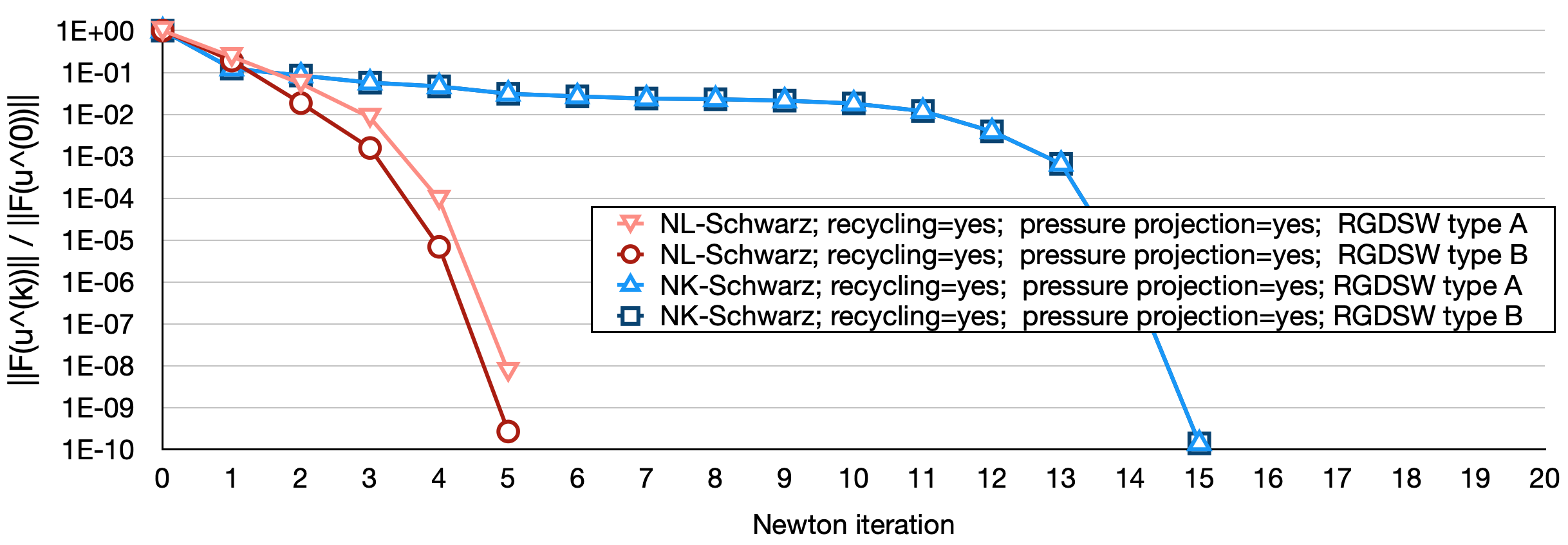}\\
\includegraphics[width=0.9\textwidth]{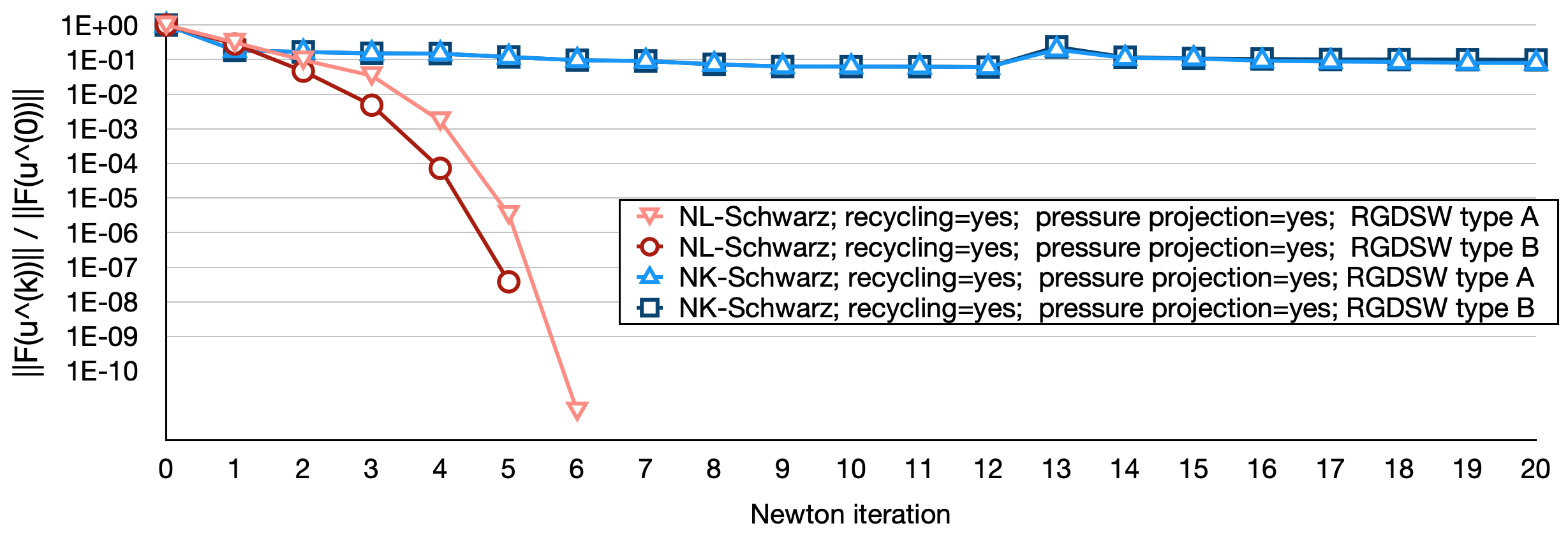}
\caption{Convergence behavior of nonlinear two-level Schwarz and Newton-Krylov-Schwarz for $Re=1000$ {\bf (top)} and $Re=1500$ {\bf (bottom)} for the lid-driven cavity flow comparing the two different RGDSW coarse spaces.}	
\label{fig:res2}
\end{figure}

To increase the robustness of Newton's method, we use a simple backtracking approach for inexact Newton methods suggested in~\cite{eisenstat}. We use this approach for the local and coarse inner loops in nonlinear Schwarz as well as for all classical Newton-Krylov-Schwarz approaches. The step size cannot get smaller then $0.01$ in the line search, which gave the best results for out specific problem. We always use 256 subdomains, since using less subdomains leads to very small coarse spaces which have no positive effect in the case of nonlinear two-level Schwarz methods. Using more subdomains is too time consuming using our MATLAB implementation. We use 128 finite elements (Taylor-Hood elements) to discretize each nonoverlapping subdomain and an overlap of $\delta=3h$. As a consequence, the original global problem has $148\,739$ degrees of freedom in total. In~\cref{table:res1} and \cref{fig:res2} we present the obtained results for the different combinations.
 
\begin{table}
\begin{center}
\caption{Comparison of hybrid restricted two-level nonlinear Schwarz {\bf (NL-Schwarz)} and hybrid restricted two-level Newton-Krylov-Schwarz method {\bf (NK-Schwarz)} with and without recycling of the coarse basis, with and without local pressure projections for the lid-driven cavity flow problem, and with two different types of the RGDSW coarse space; {\bf outer it.} gives the number of global Newton iterations; {\bf inner it.} gives the number of local Newton iterations summed up over the outer Newton iterations (minimum/maximum/average over subdomains); {\bf coarse it.} gives the number of nonlinear iterations on the second level summed up over the outer Newton iterations; {\bf GMRES it.} gives the average number of GMRES iterations per outer Newton iterations. In one case, the X marks the divergence of the method caused by the divergence of the inner coarse solve.}
\begin{tabular}{|r|r|r|r|r||c|c|c|c|}
\hline 
\multicolumn{9}{|c|}{\bf Lid-driven cavity flow} \\
\multicolumn{9}{|c|}{ $H/h=8$; {overlap} $\delta = 3h$; $N=256$ square subdomains}\\[1ex]\hline
&\bf  & \bf &  \bf pressure & \bf RGDSW & \bf outer & \bf inner it. & \bf coarse & \bf GMRES\\
$Re$ &\bf solver & \bf recycling & \bf projection & \bf type & \bf it. & \bf (min / max /avg) & \bf it. & \bf it. (avg.)\\\hline
&NL-Schwarz&no&no&A&4&9 / 11 / 9.4&14&67.0\\[0ex] 
&NL-Schwarz&yes&no&A&4&9 / 11 / 9.1&14&63.0\\[0ex] 
&NL-Schwarz&yes&yes&A&4&9 / 11 / 9.1&14&37.5\\[0ex]
500&NL-Schwarz&yes&yes&B&4&8 / 12 / 9.3&14&51.0\\[1ex]
&NK-Schwarz&yes&yes&A&7&-&-&36.4\\[0ex]
&NK-Schwarz&yes&yes&B&7&-&-&48.3\\[0ex]\hline
&NL-Schwarz&no&no&A&5&11 / 15 / 11.5&19&74.8\\[0ex] 
&NL-Schwarz&yes&no&A&5&11 / 14 / 11.5&21&69.4\\[0ex] 
&NL-Schwarz&yes&yes&A&5&11 / 14 / 11.4&20&45.2\\[0ex]
1000&NL-Schwarz&yes&yes&B&5&10 / 14 / 12.3&18&57.2\\[1ex]
&NK-Schwarz&yes&yes&A&15&-&-&44.9\\[0ex]
&NK-Schwarz&yes&yes&B&15&-&-&57.9\\[0ex]\hline
&NL-Schwarz&no&no&A&6&13 / 19 / 13.8&25&83.3\\[0ex] 
&NL-Schwarz&yes&no&A&6&13 / 19 / 13.7&26&78.2\\[0ex] 
&NL-Schwarz&yes&yes&A&6&13 / 19 / 13.7&26&52.5\\[0ex]
1500&NL-Schwarz&yes&yes&B&5&11 / 14 / 13.2&21&65.2\\[1ex]
&NK-Schwarz&yes&yes&A&$>$20&-&-&58.7\\[0ex]
&NK-Schwarz&yes&yes&B&$>$20&-&-&76.0\\[0ex]\hline
&NL-Schwarz&yes&yes&A&6&14 / 19 / 14.5&27&59.7\\[0ex]
2000&NL-Schwarz&yes&yes&B&6&13 / 18 / 15.5&26&73.2\\\hline
&NL-Schwarz&yes&yes&A&X&X&X&X\\[0ex]
2500&NL-Schwarz&yes&yes&B&7&16 / 21 / 18.2&30&79.9\\\hline
\end{tabular}
\label{table:res1}
\end{center}
\end{table}

In general, we can observe that nonlinear Schwarz methods can improve the nonlinear convergence, especially for high Reynolds numbers where the classical Newton approach does not converge within 20 iterations. For small Reynolds numbers the global Newton method is, as expected, the best choice, since it provides comparably fast convergence without the additional work introduced by the inner loops. The best setup in our tests was to always recycle the coarse basis functions and to use the local pressure projection $\overline{\mathcal{P}}_i$. The latter ones improve the linear convergence without deteriorating the fast convergence of the nonlinear solver. Comparing both RGDSW coarse spaces, the nonlinear convergence is very similar but slightly better for type B coarse spaces up to a Reynolds number of 2000. Additionally, the inner coarse loop converges a bit faster. On the other hand, the type A coarse space leads to a better linear convergence. For the case with $Re=2500$, nonlinear Schwarz with type A coarse basis functions divergence caused by a divergent inner course iteration. To summarize, type B RGDSW coarse basis functions proved to be more robust in our experiments. Additionally, the type B coarse spaces can be built without geometry information on the subdomains while for the type A coarse basis functions the coordinates of the nodes have to be known. For a more detailed picture of the convergence behavior, we provide~\cref{fig:res2} including the most promising variants for Reynolds numbers of 1000 and 1500. 

\begin{figure}
\centering
\includegraphics[width=0.3\textwidth]{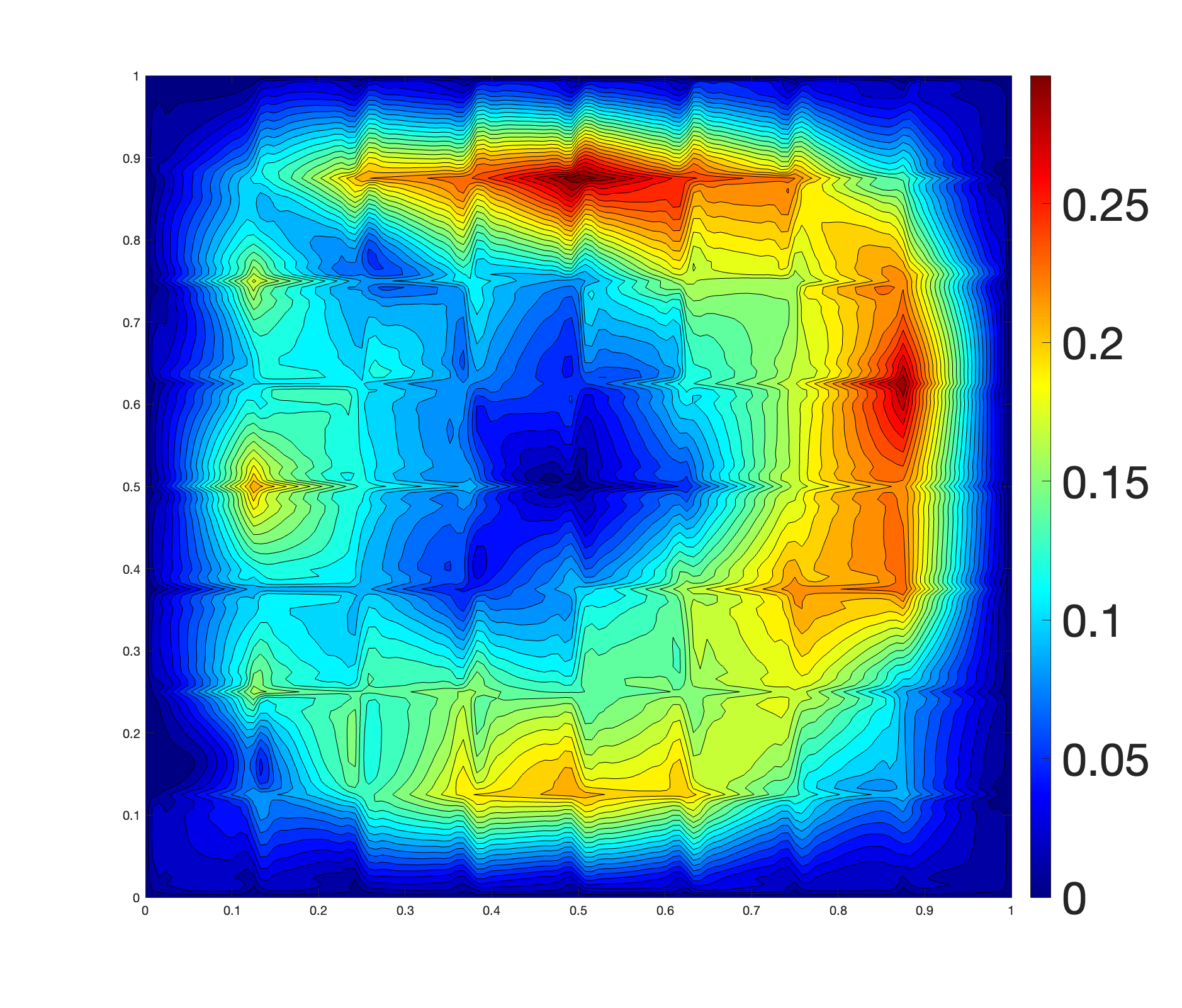}
\includegraphics[width=0.3\textwidth]{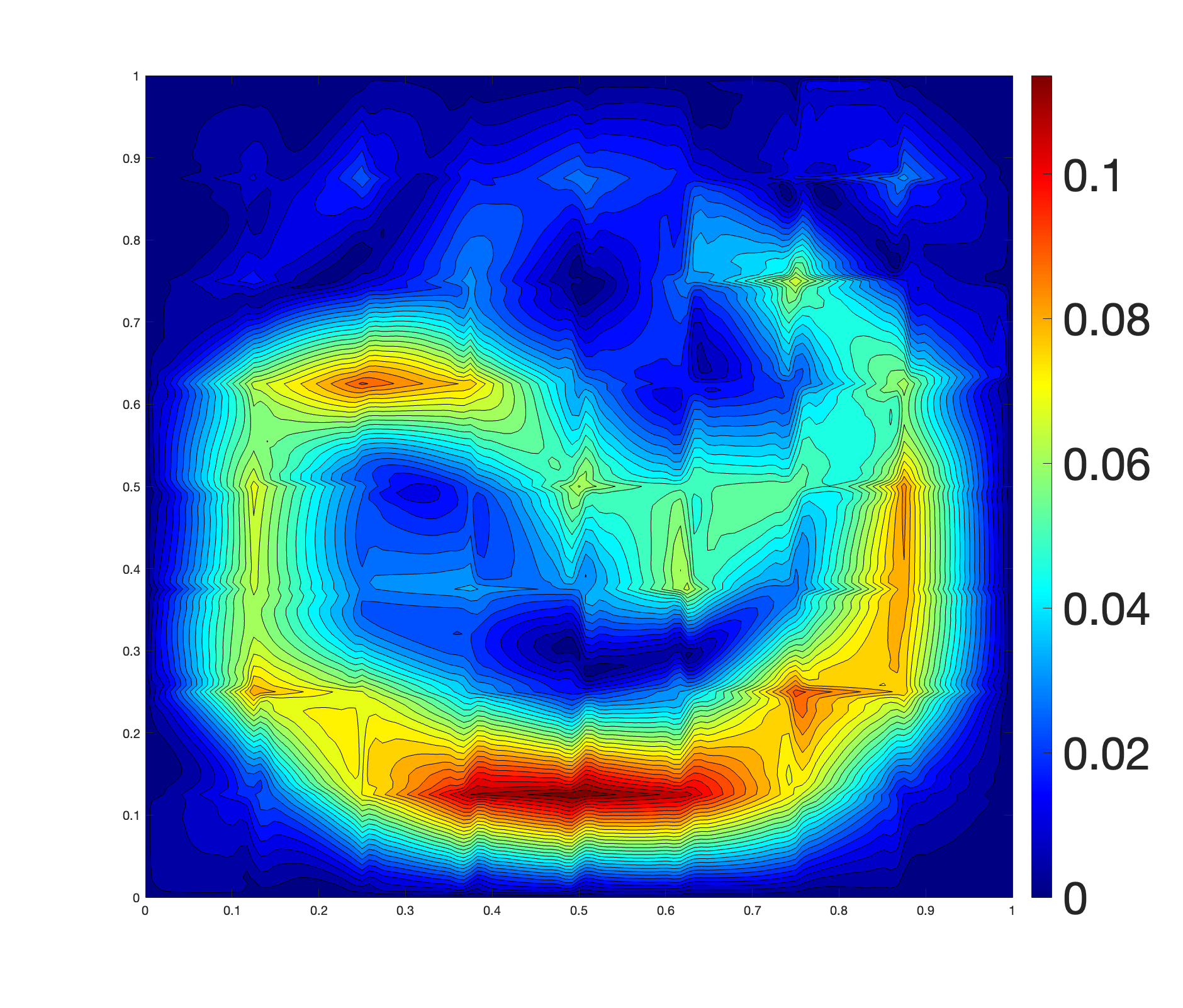}
\includegraphics[width=0.3\textwidth]{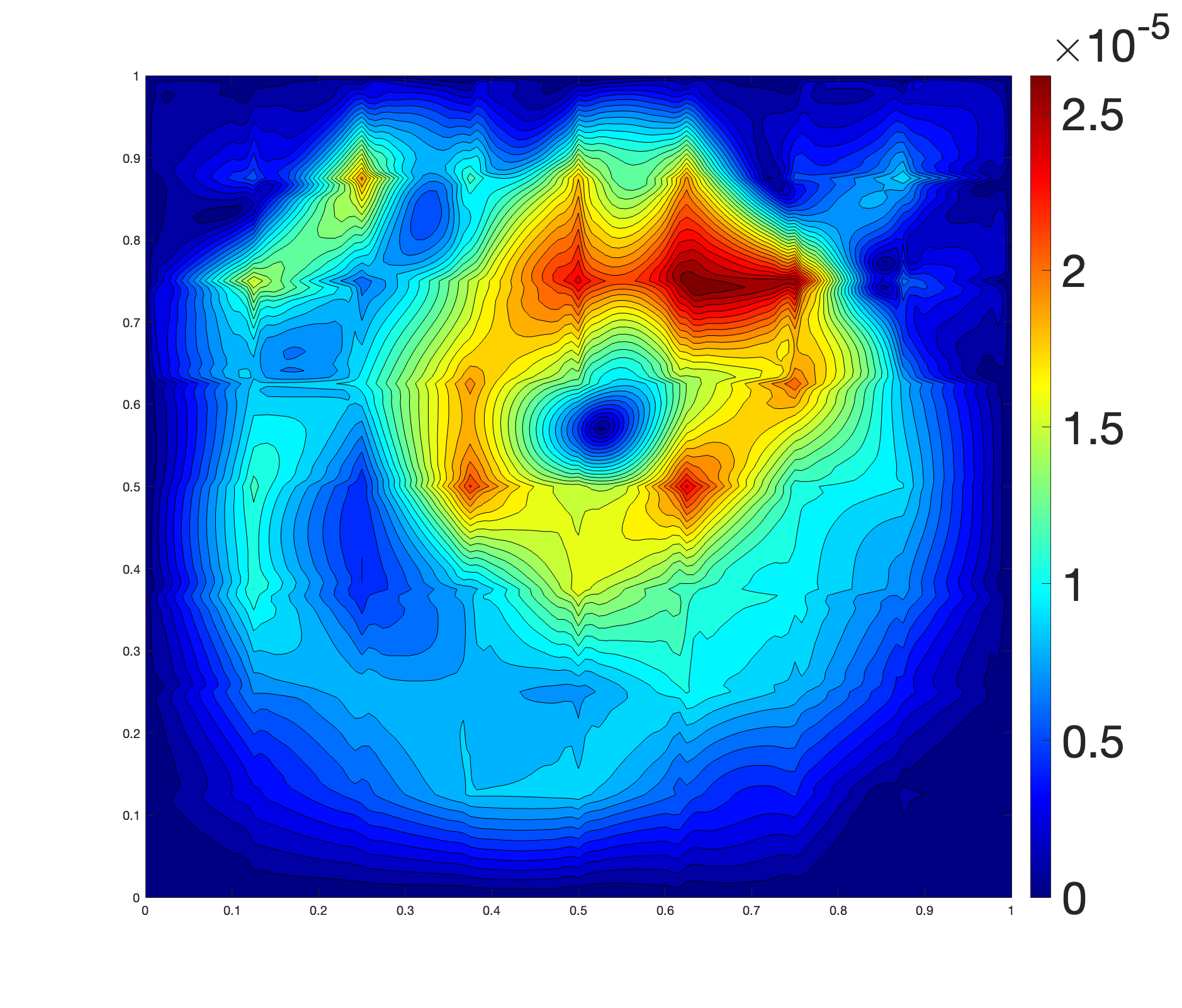}\\
\includegraphics[width=0.3\textwidth]{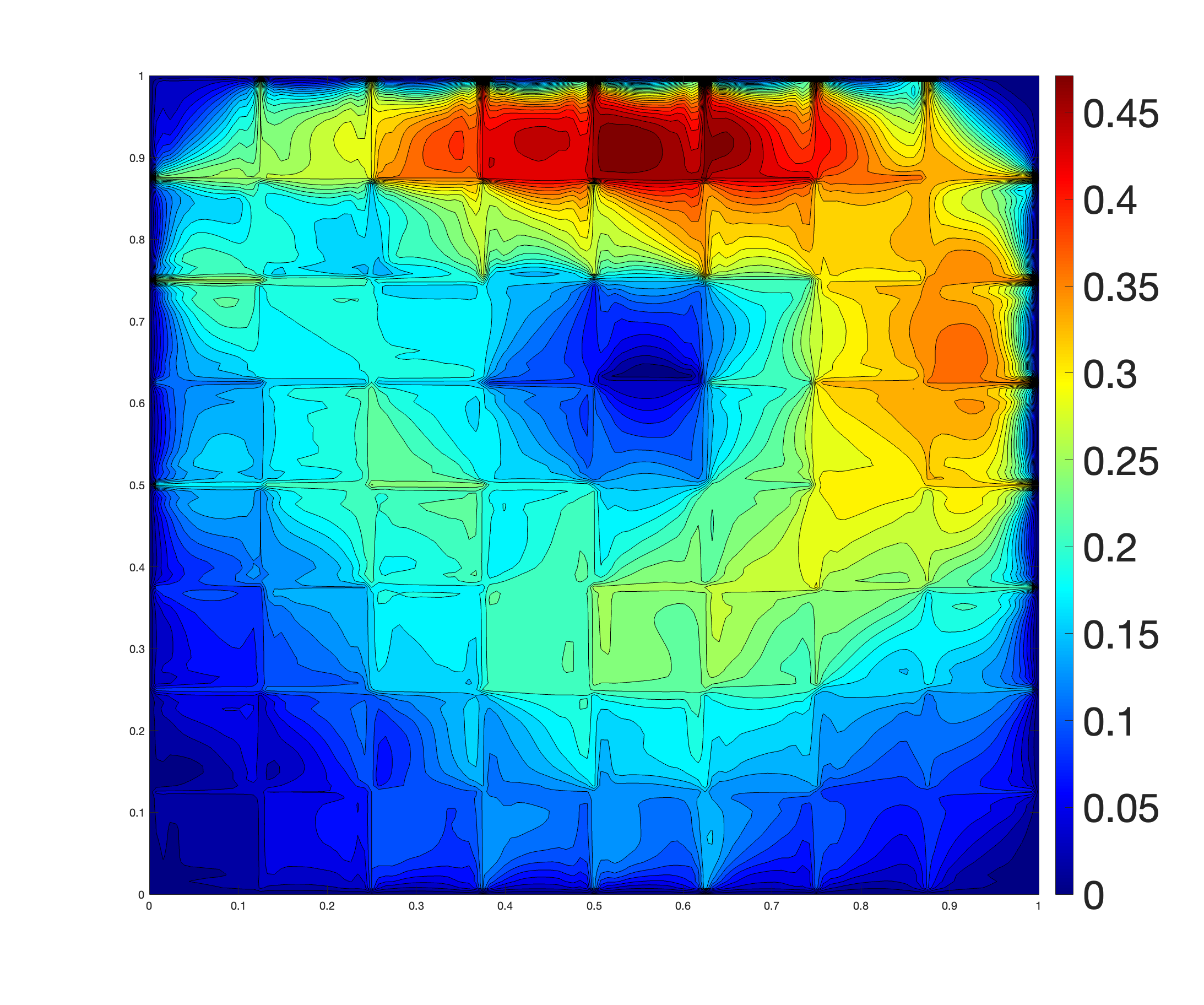}
\includegraphics[width=0.3\textwidth]{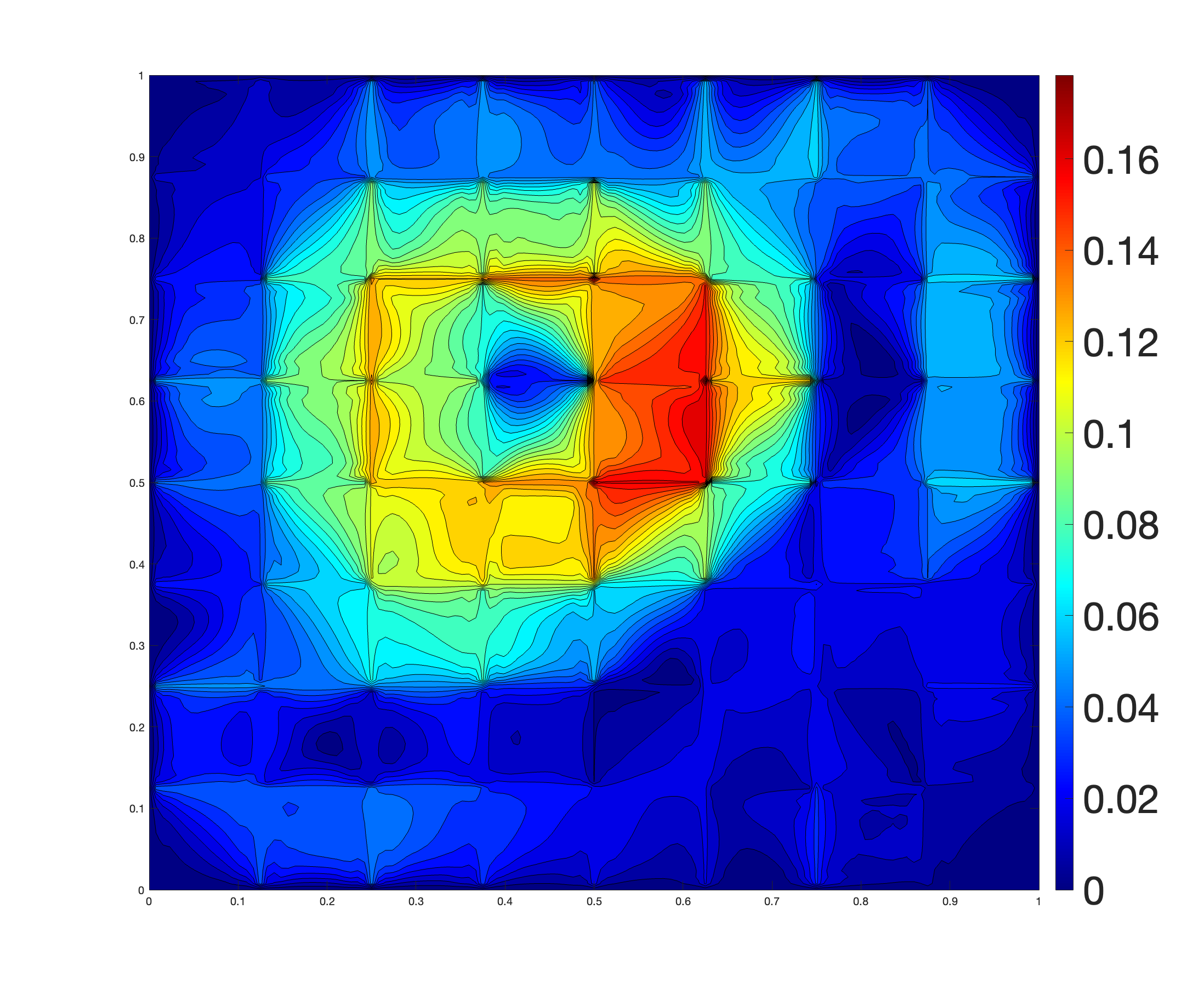}
\includegraphics[width=0.3\textwidth]{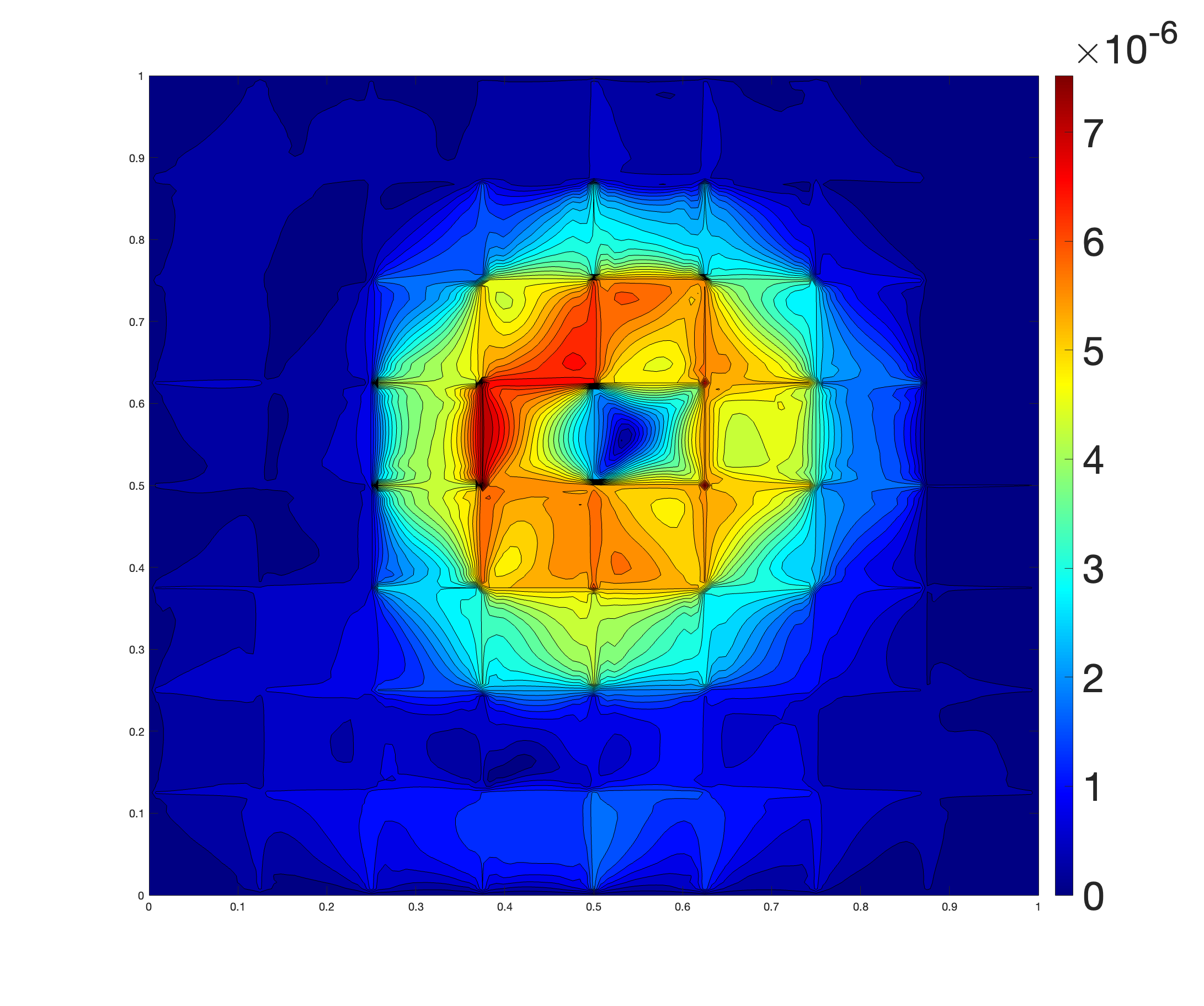}
\caption{{\bf From left to right:} Evolution of the coarse correction $T_0(u^{(k)})$; shown are the first, second, and fifth (last) iterations. The point-wise Euclidian norm of the velocity is plotted. {\bf Top row:} Type A coarse space. {\bf Bottom row:} Type B coarse space.}
\label{fig:correction}	
\end{figure}

Furthermore, we provide a visualization of the coarse corrections $T_0$ for a Reynolds number of $500$ and a decomposition into $64$ subdomains in \cref{fig:correction}. Although the convergence speed is nearly identical for both coarse spaces, it can be seen that, in general, the convergence strongly depends on the choice of the coarse basis and different coarse corrections are computed. 

\begin{figure}
\centering
\includegraphics[width=0.3\textwidth]{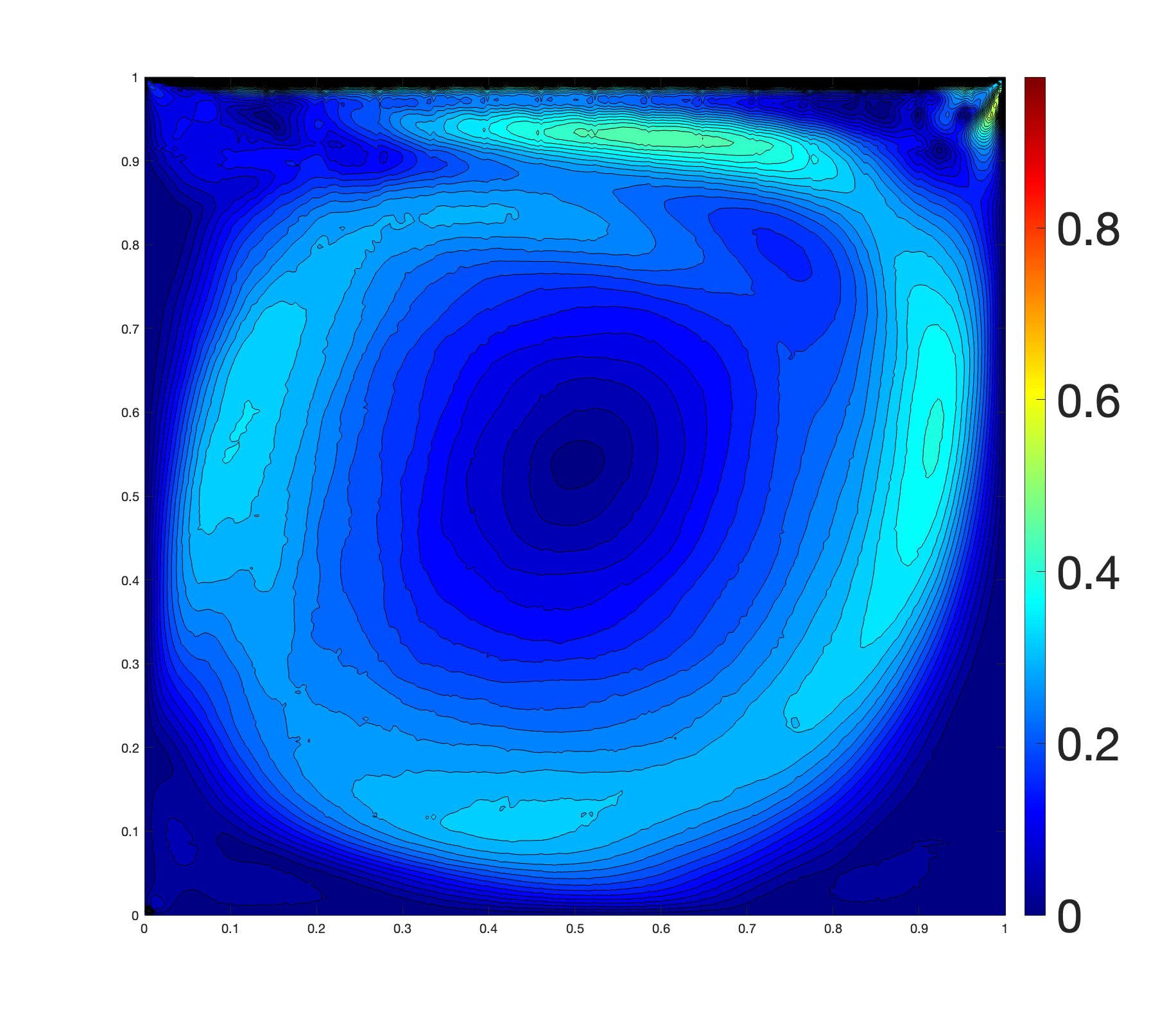}
\includegraphics[width=0.3\textwidth]{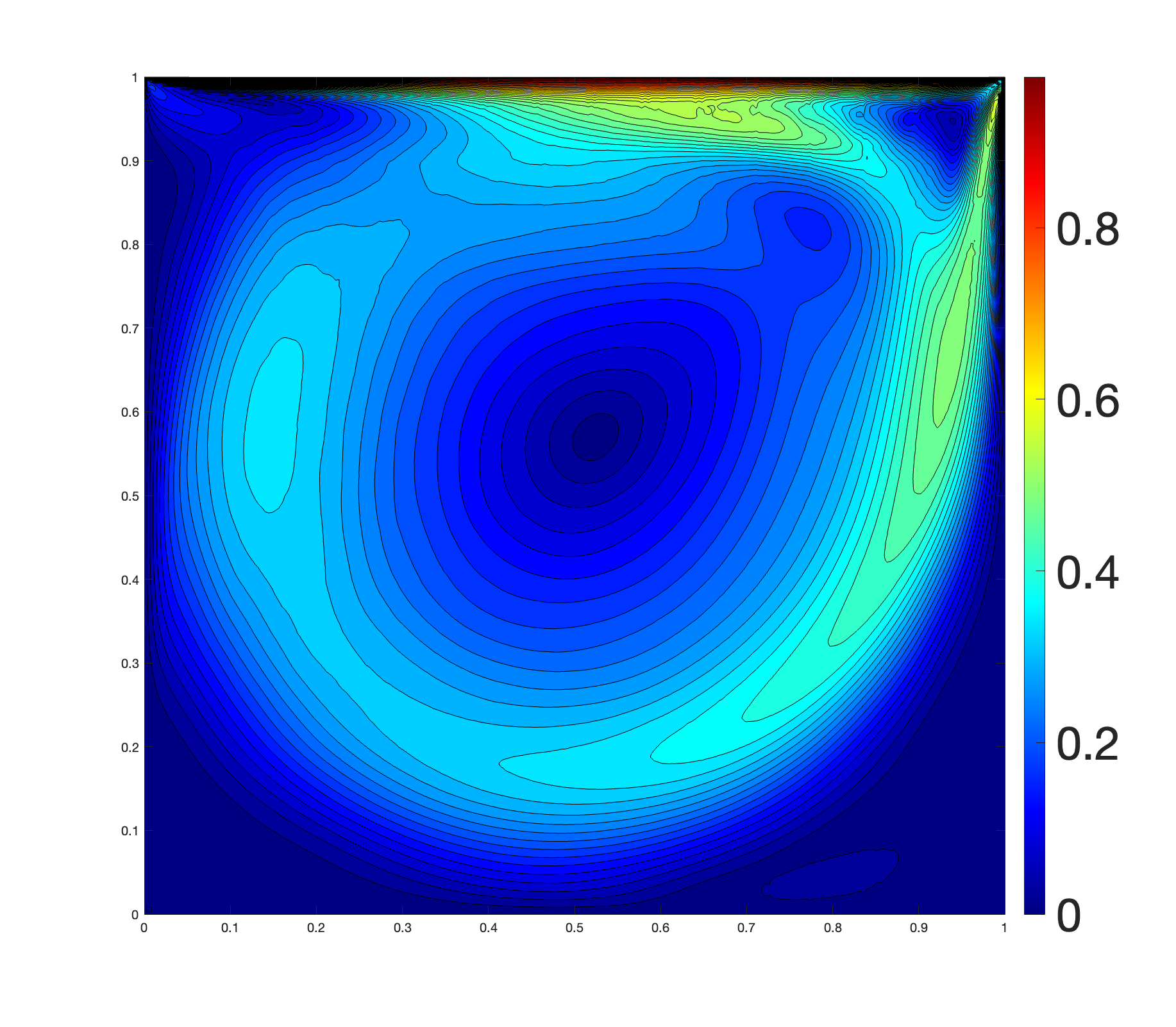}
\includegraphics[width=0.3\textwidth]{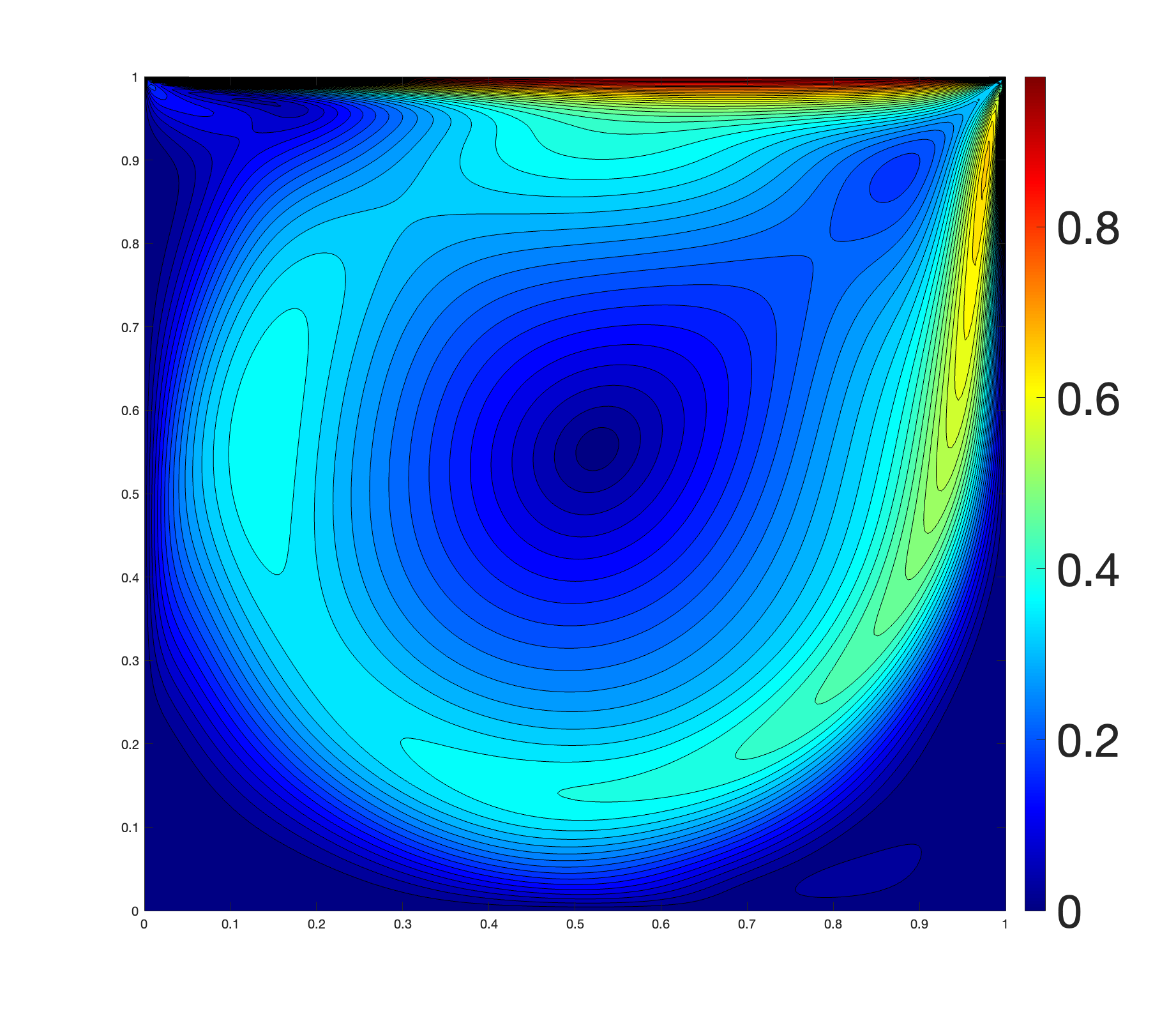}\\
\includegraphics[width=0.3\textwidth]{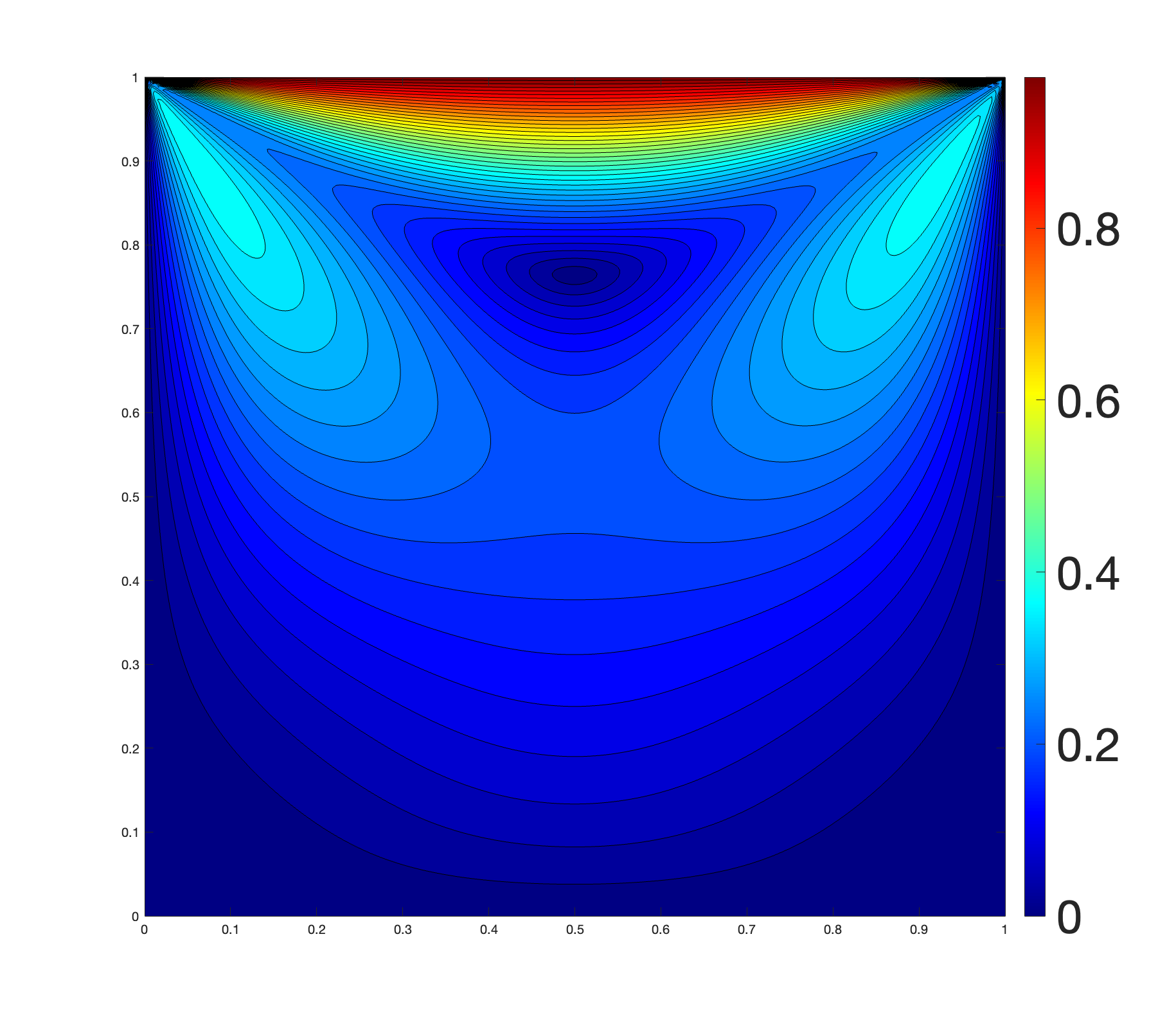}
\includegraphics[width=0.3\textwidth]{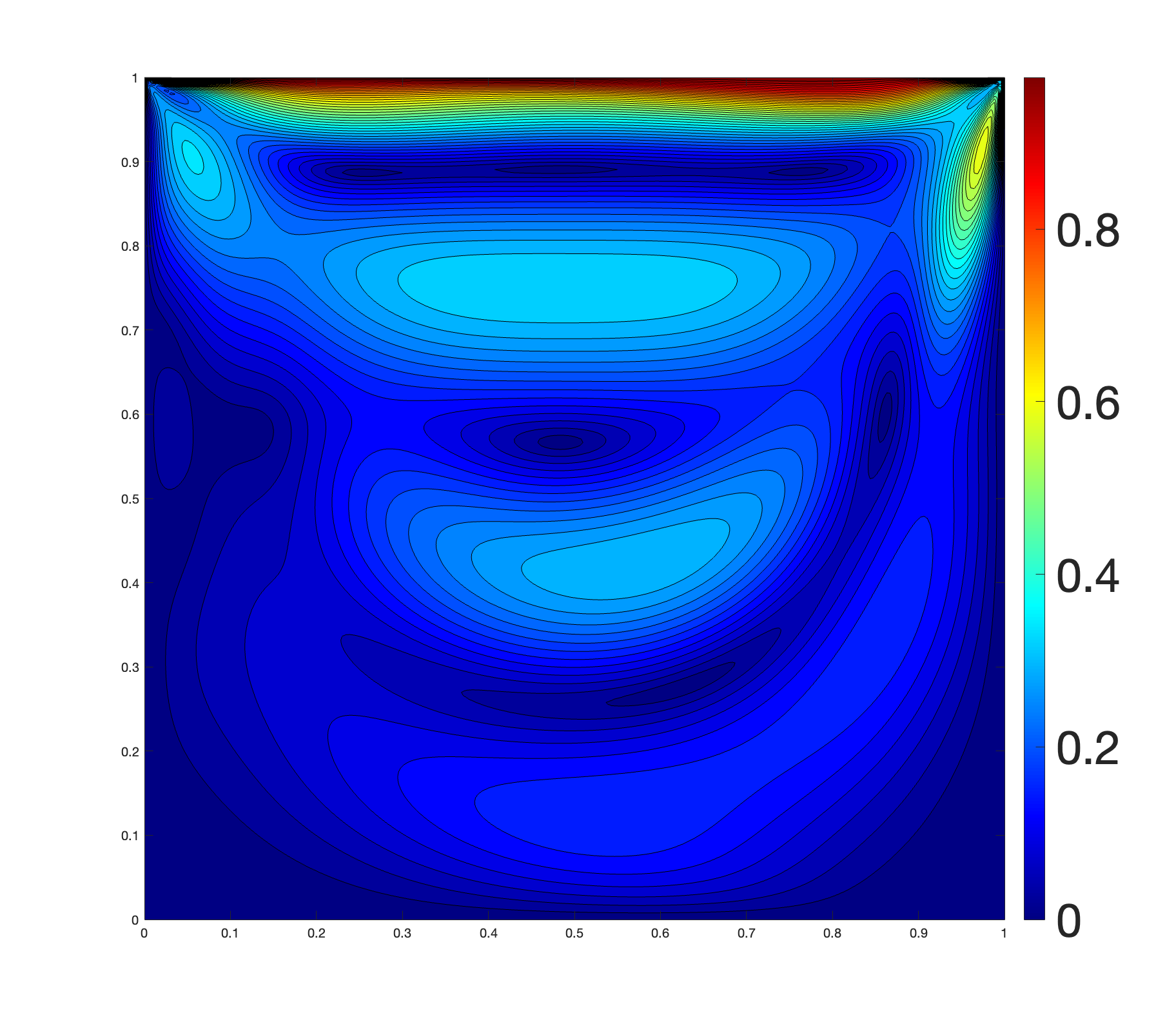}
\includegraphics[width=0.3\textwidth]{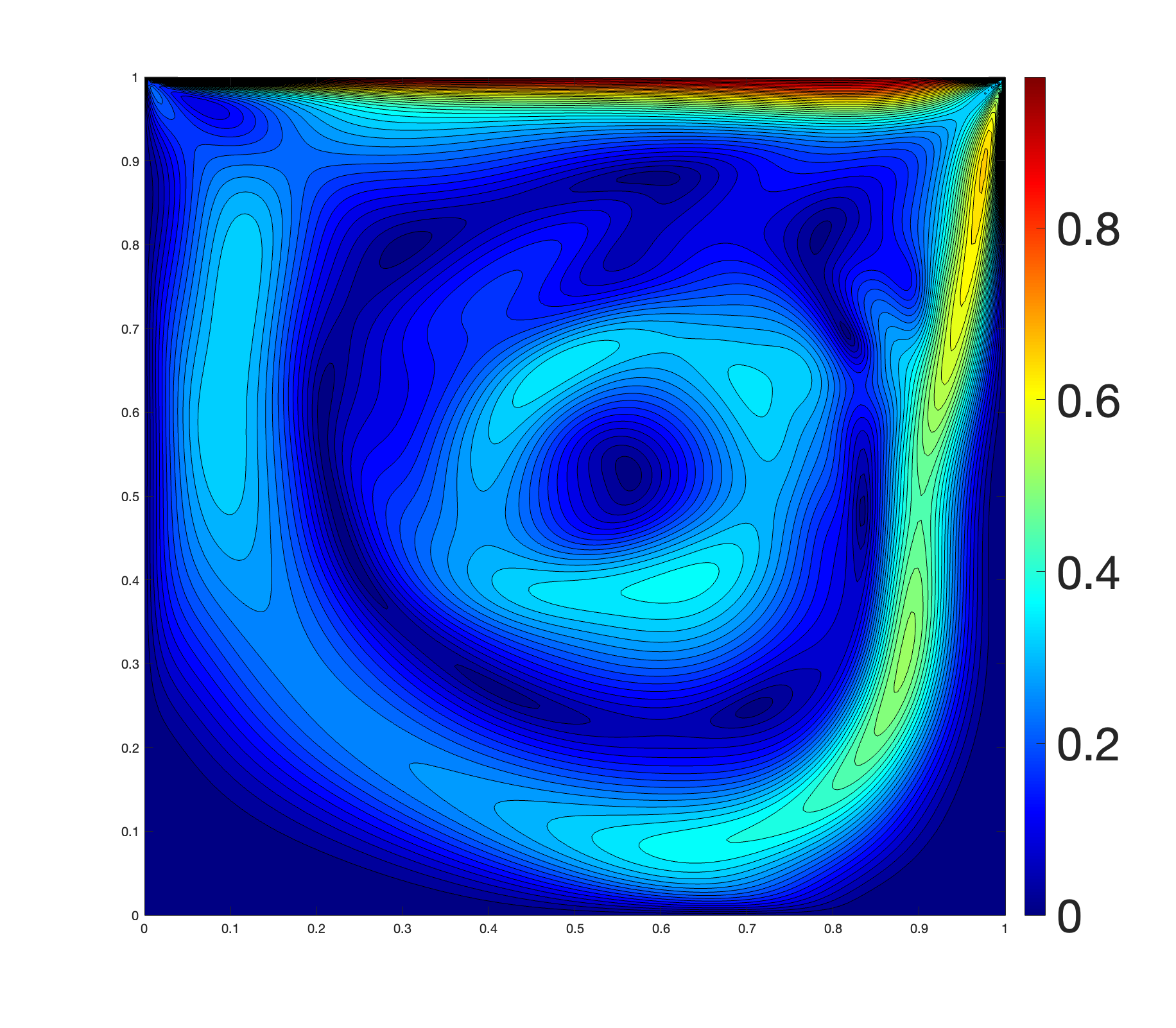}
\caption{{\bf From left to right:} Evolution of the Newton iteration; shown are the first, second, and sixth iterations for nonlinear Schwarz {\bf (top row)} and Newton-Krylov-Schwarz {\bf (bottom row)}. The point-wise Euclidian norm of the velocity is plotted.}
\label{fig:sol}	
\end{figure}

Finally, we show the Newton iterates for the case with $Re=1\,500$ for Newton-Krylov-Schwarz and nonlinear Schwarz to compare the convergence behavior; see \cref{fig:sol}. Already after one iteration, nonlinear Schwarz visibly converges against the solution while Newton-Krylov-Schwarz tends to form more vortices and clearly diverges, although a backtracking approach is used to control the step lengths.

\begin{acknowledgement}
The authors acknowledge the financial support by the German Federal Ministry of Education and Research (BMBF) for the project {\it Stroemungsraum} within the exascale computing programme SCALEXA.	
\end{acknowledgement}

\bibliographystyle{spmpsci}     
\bibliography{NL_bib}

\begin{thebibliography}{10}
\providecommand{\url}[1]{{#1}}
\providecommand{\urlprefix}{URL }
\expandafter\ifx\csname urlstyle\endcsname\relax
  \providecommand{\doi}[1]{DOI~\discretionary{}{}{}#1}\else
  \providecommand{\doi}{DOI~\discretionary{}{}{}\begingroup
  \urlstyle{rm}\Url}\fi

\bibitem{schwarz5}
Cai, X.C., Keyes, D.E.: Nonlinearly preconditioned inexact {N}ewton algorithms.
\newblock SIAM J. Sci. Comput. \textbf{24}(1), 183--200 (2002).
\newblock \doi{10.1137/S106482750037620X}.
\newblock \urlprefix\url{https://doi.org/10.1137/S106482750037620X}

\bibitem{schwarz6}
Cai, X.C., Keyes, D.E., Marcinkowski, L.: Non-linear additive {S}chwarz
  preconditioners and application in computational fluid dynamics.
\newblock pp. 1463--1470 (2002).
\newblock \doi{10.1002/fld.404}.
\newblock \urlprefix\url{https://doi.org/10.1002/fld.404}.
\newblock LMS Workshop on Domain Decomposition Methods in Fluid Mechanics
  (London, 2001)

\bibitem{gdsw1}
Dohrmann, C.R., Klawonn, A., Widlund, O.B.: Domain decomposition for less
  regular subdomains: Overlapping {S}chwarz in two dimensions.
\newblock SIAM Journal on Numerical Analysis \textbf{46}(4), 2153--2168 (2008).
\newblock \doi{10.1137/070685841}.
\newblock \urlprefix\url{https://doi.org/10.1137/070685841}

\bibitem{rgdsw}
Dohrmann, C.R., Widlund, O.B.: On the design of small coarse spaces for domain
  decomposition algorithms.
\newblock SIAM Journal on Scientific Computing \textbf{39}(4), A1466--A1488
  (2017).
\newblock \doi{10.1137/17M1114272}.
\newblock \urlprefix\url{https://doi.org/10.1137/17M1114272}

\bibitem{raspen}
Dolean, V., Gander, M.J., Kheriji, W., Kwok, F., Masson, R.: Nonlinear
  preconditioning: how to use a nonlinear {S}chwarz method to precondition
  {N}ewton's method.
\newblock SIAM J. Sci. Comput. \textbf{38}(6), A3357--A3380 (2016).
\newblock \doi{10.1137/15M102887X}.
\newblock \urlprefix\url{https://doi.org/10.1137/15M102887X}

\bibitem{eisenstat}
Eisenstat, S.C., Walker, H.F.: Globally convergent inexact {N}ewton methods.
\newblock SIAM J. Optim. \textbf{4}(2), 393--422 (1994).
\newblock \doi{10.1137/0804022}.
\newblock \urlprefix\url{https://doi.org/10.1137/0804022}

\bibitem{schwarz2}
Heinlein, A., Hochmuth, C., Klawonn, A.: Monolithic overlapping {S}chwarz
  domain decomposition methods with {GDSW} coarse spaces for incompressible
  fluid flow problems.
\newblock SIAM J. Sci. Comput. \textbf{41}(4), C291--C316 (2019).
\newblock \doi{10.1137/18M1184047}.
\newblock \urlprefix\url{https://doi.org/10.1137/18M1184047}

\bibitem{schwarz3}
Heinlein, A., Hochmuth, C., Klawonn, A.: Reduced dimension {GDSW} coarse spaces
  for monolithic {S}chwarz domain decomposition methods for incompressible
  fluid flow problems.
\newblock Internat. J. Numer. Methods Engrg. \textbf{121}(6), 1101--1119
  (2020).
\newblock \doi{10.1002/nme.6258}.
\newblock \urlprefix\url{https://doi.org/10.1002/nme.6258}

\bibitem{schwarz4}
Heinlein, A., Klawonn, A., Lanser, M.: Adaptive nonlinear domain decomposition
  methods with an application to the {$p$}-{L}aplacian.
\newblock SIAM J. Sci. Comput. \textbf{45}(3), S152--S172 (2023).
\newblock \doi{10.1137/21M1433605}.
\newblock \urlprefix\url{https://doi.org/10.1137/21M1433605}

\bibitem{schwarz1}
Heinlein, A., Lanser, M.: Additive and hybrid nonlinear two-level {S}chwarz
  methods and energy minimizing coarse spaces for unstructured grids.
\newblock SIAM J. Sci. Comput. \textbf{42}(4), A2461--A2488 (2020).
\newblock \doi{10.1137/19M1276972}.
\newblock \urlprefix\url{https://doi.org/10.1137/19M1276972}

\bibitem{nlfeti1}
Klawonn, A., Lanser, M., Rheinbach, O.: Nonlinear {FETI}-{DP} and {BDDC}
  methods.
\newblock SIAM J. Sci. Comput. \textbf{36}(2), A737--A765 (2014).
\newblock \doi{10.1137/130920563}.
\newblock \urlprefix\url{https://doi.org/10.1137/130920563}

\bibitem{nlfeti4}
K\"ohler, S., Rheinbach, O.: Globalization of nonlinear {FETI}-{DP} domain
  decomposition methods using an {SQP} approach.
\newblock Vietnam J. Math. \textbf{50}(4), 1053--1079 (2022).
\newblock \doi{10.1007/s10013-022-00567-2}.
\newblock \urlprefix\url{https://doi.org/10.1007/s10013-022-00567-2}

\bibitem{nlfeti2}
K{\"o}hler, S., Rheinbach, O.: Composing {T}wo {D}ifferent {N}onlinear
  {FETI}--{DP} {M}ethods.
\newblock In: Domain {D}ecomposition {M}ethods in {S}cience and {E}ngineering
  {XXVII}, \emph{Lect. Notes Comput. Sci. Eng.}, vol. 149, pp. 479--486.
  Springer, Cham (2024).
\newblock \doi{10.1007/978-3-031-50769-4\_57}.
\newblock \urlprefix\url{https://doi.org/10.1007/978-3-031-50769-4\_57}

\bibitem{krause}
Kopani\v{c}\'akov\'a, A., Kothari, H., Krause, R.: Nonlinear field-split
  preconditioners for solving monolithic phase-field models of brittle
  fracture.
\newblock Comput. Methods Appl. Mech. Engrg. \textbf{403}, Paper No. 115733, 31
  (2023).
\newblock \doi{10.1016/j.cma.2022.115733}.
\newblock \urlprefix\url{https://doi.org/10.1016/j.cma.2022.115733}

\bibitem{mspin1}
Liu, L., Gao, W., Yu, H., Keyes, D.E.: Overlapping multiplicative {S}chwarz
  preconditioning for linear and nonlinear systems.
\newblock J. Comput. Phys. \textbf{496}, Paper No. 112548, 22 (2024).
\newblock \doi{10.1016/j.jcp.2023.112548}.
\newblock \urlprefix\url{https://doi.org/10.1016/j.jcp.2023.112548}

\bibitem{mspin2}
Liu, L., Keyes, D.E.: Field-split preconditioned inexact {N}ewton algorithms.
\newblock SIAM J. Sci. Comput. \textbf{37}(3), A1388--A1409 (2015).
\newblock \doi{10.1137/140970379}.
\newblock \urlprefix\url{https://doi.org/10.1137/140970379}

\bibitem{nlfeti5}
Negrello, C., Gosselet, P., Rey, C.: Nonlinearly preconditioned {FETI} solver
  for substructured formulations of nonlinear problems.
\newblock Mathematics \textbf{9}(24) (2021).
\newblock \doi{10.3390/math9243165}.
\newblock \urlprefix\url{https://www.mdpi.com/2227-7390/9/24/3165}

\bibitem{nlfeti3}
Pebrel, J., Rey, C., Gosselet, P.: A nonlinear dual-domain decomposition
  method: Application to structural problems with damage.
\newblock International Journal for Multiscale Computational Engineering
  \textbf{6}(3), 251--262 (2008)

\end{thebibliography}

\end{document}